\documentclass[11pt]{article}
\usepackage{enumerate}
\usepackage{amssymb,a4wide,latexsym,makeidx,epsfig,fleqn}
\usepackage{amsthm}
\usepackage{amsmath}
\usepackage{enumerate}
\newtheorem{theorem}{Theorem}[section]

\newtheorem{lemma}[theorem]{Lemma}

\newtheorem{corollary}[theorem]{Corollary}

\begin{document}
\textwidth 150mm \textheight 225mm
\title{Upper bounds on the Q-spectral radius of book-free and/or $K_{s,t}$-free graphs
\thanks{ Supported by the National Natural Science Foundation of China (No. 11171273) and sponsored by the Seed Foundation of Innovation and Creation for Graduate Students
in Northwestern Polytechnical University (No. Z2016170).
 \vskip 0.05in}}
\author{{Qi Kong$^1$, Ligong Wang$^2$}\\
{\small Department of Applied Mathematics, School of Science, Northwestern
Polytechnical University,}\\ {\small  Xi'an, Shaanxi 710072,
People's Republic
of China.}\\ {\small $^1$E-mail: kongqixgd@163.com}\\
{\small $^2$E-mail: lgwangmath@163.com}}
\date{}
\maketitle
\begin{center}
\begin{minipage}{120mm}
\vskip 0.3cm
\begin{center}
{\small {\bf Abstract}}
\end{center}
{\small In this paper, we prove two results about the signless Laplacian spectral radius $q(G)$ of a graph $G$ of order $n$ with maximum degree $\Delta$. Let $B_{n}=K_{2}+\overline{K_{n}}$ denote a book, i.e., the graph $B_{n}$ consists of $n$ triangles sharing an edge.

(1) Let $1< k\leq l< \Delta < n$ and $G$ be a connected \{$B_{k+1},K_{2,l+1}$\}-free graph of order $n$ with maximum degree $\Delta$. Then
$$\displaystyle q(G)\leq \frac{1}{4}[3\Delta+k-2l+1+\sqrt{(3\Delta+k-2l+1)^{2}+16l(\Delta+n-1)}.$$
with equality holds if and only if $G$ is a strongly regular graph with parameters ($\Delta$, $k$, $l$).

(2) Let $s\geq t\geq 3$, and let $G$ be a connected $K_{s,t}$-free graph of order $n$ $(n\geq s+t)$. Then
$$q(G)\leq n+(s-t+1)^{1/t}n^{1-1/t}+(t-1)(n-1)^{1-3/t}+t-3.$$

 }
\end{minipage}
\end{center}

\section{Introduction }
\label{sec:ch6-introduction}

Our graph notation follows Bollob\'{a}s \cite{MGT}. In particular, let $G=(V(G), E(G))$ be a simple graph. Denote by $v(G)$ the order of $G$ and $e(G)$ the size of $G$, that is to say, $v(G)=|V(G)|$, and $e(G)=|E(G)|$. Set $\Gamma_{G}(u)=\{v|uv\in E(G)\}$, and $d_{G}(u)=|\Gamma_{G}(u)|$, or simply $\Gamma(u)$ and $d(u)$, respectively. Let $\delta = \delta(G)$ and $\Delta = \Delta(G)$ denote the minimal degree and maximal degree of graph $G$, respectively.

For a simple graph $G$ of order $n$, let the matrix $D(G)=diag(d_{1}, d_{2}, \cdots, d_{n})$, and $A(G)=(a_{ij})_{n\times n}$ be the adjacency matrix of $G$ with $a_{ij}=1$ if $v_{i}$ is adjacent to $v_{j}$, and $a_{ij}=0$ otherwise. The matrix $Q(G)=D(G)+A(G)$ is called the signless Laplacian matrix of $G$. The largest eigenvalue of $A(G)$ and $Q(G)$ are called spectral radius and signless Laplacian spectral radius (or simply Q-spectral radius) of $G$ , respectively, and marked $\rho(G)$ and $q(G)$, respectively.

Let $X$ be a set of vertices of $G$, $G[X]$ is  the graph induced by $X$, and $e(X)=e(G[X])$. Let $P_{k}$, $C_{k}$ and $K_{k}$ be the path, cycle, and complete graph of order $k$, respectively. If all vertices of $G$ have the same degree $k$, then $G$ is \textit{$k$-regular}. A $k$-regular graph is called \textit{strongly regular} with parameters $(k, a, c)$ whenever each pair of adjacent vertices have $a\geq 0$ common neighbors, and each pair of non-adjacent vertices have $c\geq 1$ common neighbors.

The main results of this paper are in the spirt of the trend in the famous Zarankiewicz problem \cite{Snr}:

\textbf{Problem A} \textit{How many edges can have a graph of order $n$ if it does not contain a complete bipartite subgraph $K_{s, t}$ ?}

In 1996, F\"{u}redi \cite{Aub} gave an upper bound on the above Zarankiewicz problem. In 2010, Nikiforov \cite{Act} improved his result. That is, if $G$ is a $K_{s, t}$ -free graph of order $n$, then
$$e(G)\leq \frac{1}{2}(s-t+1)^{1/t}n^{2-1/t}+\frac{1}{2}(t-1)n^{2-2/t}+\frac{1}{2}(t-2)n.$$
The spectral version of the Zarankiewicz problem is the following one:

\textbf{Problem B} \textit{How large can be the spectral radius $\rho(G)$ of a graph $G$ of order $n$ that does not contain $K_{s, t}$ ?}

There are some results for some value of $s$ and $t$.

In 2007, the upper bound on the signless Laplacian spectral radius of $K_{2,l+1}$-free graph as the corollary of the following Lemma \ref{le:l1} was proved in \cite{Ubo} by Shi and Song.

\begin{lemma}\label{le:l1}
$0\leq k\leq l\leq \Delta < n$ and $G$ be a connected \{$B_{k+1},K_{2,l+1}$\}-free graph of order $n$ with maximum degree $\Delta$. Then
$$\rho(G)\leq [k-l+\sqrt{(k-l)^{2}+4\Delta+4l(n-l)}]/2,$$
with equality if and only if $G$ is a stongly regular with parameters $(\Delta, k, l)$.
\end{lemma}

In 2007, Nikiforov \cite{Bog} improved the above bound showing that
\begin{lemma}\label{le:l2}
Let $l\geq k\geq 0$. If $G$ is a \{$B_{k+1},K_{2,l+1}$\}-free graph of order $n$ with maximum degree $\Delta$. Then
$$\rho(G)\leq min\{\Delta, \frac{1}{2}[k-1+1+\sqrt{(k-l+1)^{2}+4l(n-1)}]\}.$$
If $G$ is connected, equality holds if and only if one of the following conditions holds:

(1) $\Delta^{2}-\Delta(k-l+1)\leq l(n-1)$ and $G$ is $\Delta$-regular;

(2) $\Delta^{2}-\Delta(k-l+1)> l(n-1)$ and every two vertices of $G$ have $k$ common neighbors if they are adjacent, and $l$ common neighbors otherwise.
\end{lemma}

Setting $l=\Delta$ or $k=l$, Lemma \ref{le:l2} implies assertions that strengthen Corollaries 1 and 2 of  \cite{Ubo}.

In 2010, Nikiforov \cite{Act} also gave a bound as the following lemma.
\begin{lemma}\label{le:l3}
let $s\geq t\geq 2$, and let $G$ be a $K_{s, t}$-free graph of order $n$. If $t=2$, then
$$\rho(G)\leq 1/2+\sqrt{(s-1)(n-1)+1/4}.$$
If $t\geq 3$, then
$$\rho(G)\leq(s-t+1)^{1/t}n^{1-1/t}+(t-1)n^{1-2/t}+t-2.$$
and
$$e(G)<\frac{1}{2}(s-t+1)^{1/t}n^{2-1/t}+\frac{1}{2}(t-1)n^{2-2/t}+\frac{1}{2}(t-2)n.$$
\end{lemma}

A newer trend in extremal graph theory is the Zarankiewicz problem for signless Laplacian spectral radius of graphs:

\textbf{Problem C} \textit{How large can be the signless Laplacian spectral radius $q(G)$ of a graph $G$ of order $n$ that does not contain subgraph $K_{s, t}$?}

When $s=t=2$, we notice that the $K_{2,2}$-free graph is the same as $C_{4}$-free graph. Also in 2013, de Freitus \cite{Mot} has proved that if $G$ contains no $C_{4}$, then
$$q(G)<q(F_{n}),$$
unless $G=F_{n}$, where $F_{n}$ is the friendship graph of order $n$.  For $n$ odd, $F_{n}$ is a union of $\lfloor n/2\rfloor$ triangles sharing a single common vertex, and for $n$ even, $F_{n}$ is obtained by hanging an edge to the common vertex of $F_{n-1}$.

In this paper, we discuss upper bounds on the signless Laplacian spectral radius of Book-free and/or $K_{2,l+1}$-free $(l>1)$ graphs of order $n$ with maximum degree $\Delta$.
\begin{theorem}\label{th:t1}
Let $1< k\leq l< \Delta < n$ and $G$ be a connected \{$B_{k+1},K_{2,l+1}$\}-free graph of order $n$ with maximum degree $\Delta$. Then
$$\displaystyle q(G)\leq \frac{1}{4}[3\Delta+k-2l+1+\sqrt{(3\Delta+k-2l+1)^{2}+16l(\Delta+n-1)}.\hspace{60pt}(1)$$
with equality holds if and only if $G$ is a strongly regular graph with parameters ($\Delta$, $k$, $l$).
\end{theorem}

Because every graph is obviously $K_{2,\Delta+1}$-free, Theorem \ref{th:t1} readily implies a sharp upper bound for book-free graph.

\begin{corollary}\label{co:c1}
Let $1< k<\Delta< n$ and $G$ be a connected $B_{k+1}$-free graph of order $n$ with maximum degree $\Delta$. Then
$$\displaystyle q(G)\leq \frac{1}{4}[\Delta+k+1+\sqrt{(\Delta+k+1)^{2}+32\Delta(n-1)}.$$
with equality if and only if $G$ is a strongly regular graph with parameters$ (\Delta, k, \Delta)$.
\end{corollary}

Becase a $K_{2,l}$-free graph is also $B_{l}$-free. Theorem \ref{th:t1} with $k=l$ also implies a sharp upper bound for $K_{2,l}$-free graphs.

\begin{corollary}\label{co:c2}
Let $1<l<\Delta$ and $G$ be a connected $K_{2,l+1}$-free graph of order $n$ with maximum degree $\Delta$. Then
$$\displaystyle q(G)\leq \frac{1}{4}[3\Delta-l+1+\sqrt{(3\Delta-l+1)^{2}+32l(n-1)}.$$
with equality if and only if $G$ is a strongly regular graph with parameters $(\Delta, l, l)$.
\end{corollary}

Furthermore we will discuss $s\geq t\geq 3$, let $G$ be a connected graph of order $n$, when $n< s+t$, then $G$ must contain no $K_{s,t}$, so we only discuss $n\geq s+t$.
\begin{theorem}\label{th:t2}
Let $s\geq t\geq 3$, and let $G$ be a connected $K_{s,t}$-free graph of order $n$ $(n\geq s+t)$. Then
$$q(G)\leq n+(s-t+1)^{1/t}n^{1-1/t}+(t-1)(n-1)^{1-3/t}+t-3.$$
\end{theorem}

\section{Main Lemmas}
In this section, we state some well-know results which will be used in this paper.

\begin{lemma}\label{le:c4}
Let $s\geq 2$, $t\geq 2$, $0\leq k\leq s-2$, and let $G(A, B)$ be a bipartite graph with parts $A$ and $B$. Suppose that $G$ contains no copy of $K_{s,t}$ with a vertex class of size $s$ in $A$ and a vertex class of size $t$ in $B$. Then $G(A,B)$ has at most
$$(s-k-1)^{1/t}|B||A|^{1-1/t}+(t-1)|A|^{1+k/t}+k|B|$$
edges.
\end{lemma}

\begin{lemma}\label{le:c5} (\cite{FeYu}, \cite{Mer}) For every graph $G$, we have
$$\displaystyle q(G)\leq\max\limits_{u\in V(G)}\{d(u)+\frac{1}{d(u)}\sum\limits_{v\in \Gamma(u)}d(v)\}.$$
\end{lemma}

\section{Proofs}
{\bf Proof of Theorem~\ref{th:t1}.}
 Let $Q_{i}$ denote the $i$th row vector of $Q(=Q(G))$ and let $\textbf{x}=(x_{1}, x_{2}, \ldots, x_{n})^{T}$ be the Perron-eigenvector of $Q$ corresponding to $q(G)$. Then $x_{i}>0$ for $1\leq i\leq n$. Since $G$ is \{$B_{k+1}$, $K_{2,l+1}$\}-free, each pair of adjacent vertices has at most $k$ common neighbors and each pair of non-adjacent vertices has at most $l$ common neighbors. Thus
$$\sum\limits_{i=1}\limits^{n} \sum\limits_{v_{p},v_{q}\in \Gamma(v_{i})}x_{p}x_{q}\leq k
\sum\limits_{v_{p}v_{q}\in E(G)} x_{p}x_{q}+ l\sum\limits_{v_{p}v_{q}\notin E(G)} x_{p}x_{q}.\hspace{100pt}(2)$$
Then by virtue of $\textbf{x}^{T} A(K_{n}) \textbf{x}\leq \rho(K_{n})=n-1$. Thus
\begin{align*}
q(G)&=\textbf{x}^{T} Q \textbf{x}=\textbf{x}^{T} D \textbf{x}+\textbf{x}^{T} A
\textbf{x}=\sum\limits_{i=1}\limits^{n} d_{i}x_{i}^{2}+2\sum\limits_{v_{i}v_{p}\in E(G)}x_{i}x_{p}\\
&\leq \Delta+\textbf{x}^{T} A(K_{n}) \textbf{x}-2\sum\limits_{v_{i}v_{p}\notin E(G)}x_{i}x_{p}\\
&\leq \Delta+n-1-2\sum\limits_{v_{i}v_{p}\notin E(G)}x_{i}x_{p}.
\end{align*}
Also we can obtain
\begin{align*}
q(G)&=\textbf{x}^{T} Q \textbf{x}=\sum\limits_{i=1}\limits^{n}\sum\limits_{j=1,i<j}\limits^{n}2q_{i,j}x_{i}x_{j}+\sum\limits_{i=1}\limits^{n} d_{i}x_{i}^{2}\\
&\leq\sum\limits_{i=1}\limits^{n}\sum\limits_{j=1,i<j}\limits^{n}q_{i,j}(x_{i}^{2}+x_{j}^{2})+\sum\limits_{i=1}\limits^{n}d_{i}x_{i}^{2}\\
&=\sum\limits_{i=1}\limits^{n}\sum\limits_{j=1,i<j}\limits^{n}q_{i,j}x_{i}^{2}+\sum\limits_{i=1}\limits^{n}d_{i}x_{i}^{2}\\
&=2\sum\limits_{i=1}\limits^{n}d_{i}x_{i}^{2}.
\end{align*}
So
$$\sum\limits_{i=1}\limits^{n}d_{i}x_{i}^{2}\geq \frac{q}{2}.$$
Then
\begin{align*}
q^{2}(G)&=\|Q\textbf{x}\|^{2}=\sum\limits_{i=1}\limits^{n}(Q_{i}\textbf{x})^{2}=\sum\limits_{i=1}\limits^{n}(d_{i}x_{i}+\sum\limits_{v_{i}v_{p}\in E(G)}x_{p})^{2}\\
        &=\sum\limits_{i=1}\limits^{n}[d_{i}^{2}x_{i}^{2}+2d_{i}x_{i}\sum\limits_{v_{i}v_{p}\in E(G)}x_{p}+(\sum\limits_{v_{i}v_{p}\in E(G)}x_{p})^{2}]\\
        &=\sum\limits_{i=1}\limits^{n}d_{i}^{2}x_{i}^{2}+2\sum\limits_{i=1}\limits^{n}d_{i}\sum\limits_{v_{i}v_{p}\in E(G)}x_{i}x_{p}+\sum\limits_{i=1}\limits^{n}d_{i}x_{i}^{2}+2\sum\limits_{i=1}\limits^{n}\sum\limits_{v_{p},v_{q}\in \Gamma(v_{i})}x_{p}x_{q}\\
        &\leq (\Delta+1)\sum\limits_{i=1}\limits^{n}d_{i}x_{i}^{2}+2\Delta\sum\limits_{i=1}\limits^{n}\sum\limits_{v_{i}v_{p}\in E(G)}x_{i}x_{p}\\
        &~~~~+2k\sum\limits_{v_{p}v_{q}\in E(G)}x_{p}x_{q}+2l\sum\limits_{v_{p}v_{q}\notin E(G)}x_{p}x_{q}\hspace{120pt}(3)\\
        &=(\Delta+1)\sum\limits_{i=1}\limits^{n}d_{i}x_{i}^{2}+(4\Delta+2k)\sum\limits_{v_{i}v_{p}\in E(G)}x_{i}x_{p}+2l\sum\limits_{v_{p}v_{q}\notin E(G)}x_{p}x_{q}\\
        &\leq(2\Delta+k)(\sum\limits_{i=1}\limits^{n}d_{i}x_{i}^{2}+2\sum\limits_{v_{i}v_{p}\in E(G)}x_{i}x_{p})\\
        &~~~~-(\Delta+k-1)\sum\limits_{i=1}\limits^{n}d_{i}x_{i}^{2}+2l\sum\limits_{v_{p}v_{q}\notin E(G)}x_{p}x_{q}\\
        &\leq(2\Delta+k)q-\frac{\Delta+k-1}{2}q+l(\Delta+n-1-q)\\
        &=\frac{1}{2}(3\Delta+k-2l+1)q+l(\Delta+n-1).
\end{align*}
Solving the inequality gives the upper bound
$$\displaystyle q(G)\leq \frac{1}{4}[3\Delta+k-2l+1+\sqrt{(3\Delta+k-2l+1)^{2}+16l(\Delta+n-1)}.$$
If the upper bound of (1) is attained then all inequalities in the above argument must be equalities. In particular, from (2) and $x_{i}>0$ for $1\leq i\leq n$, we have that each pair of adjacent vertices in $G$ has exactly $k$ common neighbors and each pair of non-adjacent vertices in $G$ has exactly $l$ common neighbors. Moreover, by (3), $G$ must be $\Delta$-regular. Thus $G$ must be a strongly regular graph with parameters $(\Delta, k, l)$.  $\square$

{\bf Proof of Theorem~\ref{th:t2}.} By Lemma \ref{le:c5}, let $w$ be a vertex of $G$ such that
$$\displaystyle d(w)+\frac{1}{d(w)}\sum\limits_{i\in \Gamma(w)}d(i)=\max\limits_{u\in V(G)}\{d(u)+\frac{1}{d(u)}\sum\limits_{v\in \Gamma(u)}d(v)\}.$$
Then
$$\displaystyle q(G)\leq d(w)+\frac{1}{d(w)}\sum\limits_{i\in \Gamma(w)}d(i).$$
Note first that if $d(w)\leq s+t-1$, then
\begin{align*}
 q(G)&\leq d(w)+\frac{1}{d(w)}\sum\limits_{i\in \Gamma(w)}d(i)\leq d(w)+\Delta(G)\\
     &\leq s+t-1+n-1=s+t+n-2\\
     &\leq n+(s-t+1)^{1/t}n^{1-1/t}+(t-1)(n-1)^{1-3/t}+t-3.
\end{align*}
Therefore we shall assume that $s+t-1\leq d(w)\leq n-1$. Let $U$ and $W$ be disjoint sets satisfying $|U|=d(w)$ and $|W|=n-1$, and let $\varphi_{U}$ and $\varphi_{W}$ be bijections
$$\varphi_{U}: U\rightarrow \Gamma(w), \varphi_{W}: W\rightarrow V(G)\backslash \{w\}.$$
Define a bipartite graph $H$ with vertex classes $U$ and $W$ by joining $u\in U$ and $v\in W$ whenever $\{\varphi_{U}(u),\varphi_{W}(v)\}\in E(G)$.

Then we can get that $H$ does not contain a copy of $K_{s-1, t}$ with $s-1$ vertices in $W$ and $t$ vertices in $U$. Indeed, the map $\psi: V(H)\rightarrow V(G)$ defined as
 \begin{displaymath}
\psi(x)=\left\{\
        \begin{array}{ll}
          \varphi_{U}(x),~~~~~~&\mbox {if}~x\in U,\\
          \varphi_{W}(x),~~~~~~&\mbox {if}~x\in W.\\
        \end{array}
      \right.
\end{displaymath}
is a homomorphism of $H$ into $G-w$. Assume for a contradiction that $F\subset H$ is a copy of $K_{s-1, t}$ with a set of $S$ of $s-1$ vertices in $W$ and a set of $T$ of $t$ vertices in $U$. Clearly $S$ and $T$ are the vertex classes of $F$. Note that $\psi(F)$ is a copy of $K_{s-1,t}$ in $G-w$, and $\psi(S)=\varphi_{W}(S)\subset V(G)\setminus\{w\}$ and $\psi(T)=\varphi_{U}(T)\subset \Gamma_{G}(w)$ are the vertex classes of $\psi(F)$ of size $s-1$ and size $t$, respectively. Now, adding $w$ to $\psi(F)$, we see that $G$ contains a $K_{s,t}$, a contradiction proving the clain.

Suppose that $0\leq k\leq min\{s,t\}-2$. Setting $k'=k-1, s'=s-1, t'=t, A=W, B=U$, then from Lemma \ref{le:c4}, we have
$$e(H)\leq (s-k-1)^{1/t}|U||W|^{1-1/t}+(k-1)|U|+(t-1)|W|^{1+(k-1)/t}$$
$$~~~~=(s-k-1)^{1/t}d(w)n^{1-1/t}+(k-1)d(w)+(t-1)(n-1)^{1+(k-1)/t}.$$
On the other hand. We have
$$e(H)=\sum\limits_{v\in\Gamma(w)}d(v)-d(w),$$
and so,
$$\sum\limits_{v\in\Gamma(w)}d(v)\leq((s-k-1)^{1/t}n^{1-1/t}+k)d(w)+(t-1)(n-1)^{1+(k-1)/t}.$$
And then from Lemma \ref{le:c5}, we have
\begin{align*}
q(G)&\leq d(w)+\frac{1}{d(w)}\sum\limits_{i\in \Gamma(w)}d(i)\\
&\leq d(w)+\frac{(t-1)(n-1)^{1+(k-1)/t}}{d(w)}+(s-k-1)^{1/t}n^{1-1/t}+k.
\end{align*}
Since the function
$$\displaystyle f(x)=x+\frac{(t-1)(n-1)^{1+(k-1)/t}}{x}$$
is convex for $x>0$, its maximum in any closed interval is attained at one of the ends of this interval. In the case $s+t-1\leq d(w)\leq n-1$, then,
\begin{align*}
q(G)&\leq d(w)+\frac{1}{d(w)}\sum\limits_{i\in \Gamma(w)}d(i)\\
&\leq max\{s+t-1+\frac{(t-1)(n-1)^{1+(k-1)/t}}{s+t-1}, n-1+\frac{(t-1)(n-1)^{1+(k-1)/t}}{n-1}\}\\
      &~~~~+(s-k-1)^{1/t}n^{1-1/t}+k\\
&\leq (s-k-1)^{1/t}n^{1-1/t}+k+\frac{(t-1)(n-1)^{1+(k-1)/t}}{n-1}+n-1\\
&=(s-k-1)^{1/t}n^{1-1/t}+k+(t-1)(n-1)^{(k-1)/t}+n-1.
\end{align*}
Now,if $s\geq t\geq 3$, setting $k=t-2$, we obtain
$$q(G)\leq n+(s-t+1)^{1/t}n^{1-1/t}+(t-1)(n-1)^{1-3/t}+t-3.$$
So, the proof is complete. $\square$

\end{document}